\documentclass[12pt]{article}

\begin{document}

\title{Generalized Lagrange-Weyl structures and compatible connections}
\author{Mircea Cr\^a\c sm\u areanu\thanks{Partially supported by CNCSIS Grant 2005}}
\date{}
\maketitle

\begin{abstract}
Generalized Lagrange-Weyl structures and compatible connections are
introduced as a natural generalization of similar notions from Riemannian
geometry. Exactly as in Riemannian case, the compatible connection is unique
if certain symmetry conditions with respect to vertical and horizontal
Christoffel symbols are imposed.
\end{abstract}

{\bf 2000 Math. Subject Classification}: 53C60.

{\bf Key words}: generalized Lagrange-Weyl space; compatible connection.

\section*{Introduction}

Soon after the creation of general theory of relativity, Hermann Weyl
attempted in \cite{w:h} an unification of gravitation and electromagnetism
in a model of space-time geometry combining conformal and projective
structures.

Let ${\cal G}$ be a conformal structure on the smooth manifold $M$ i.e. an
equivalence class of Riemannian metrics: $g\sim ~\overline{g}$ if there
exists a smooth function $f\in C^{\infty }\left( M\right) $ such that $%
\overline{g}=e^{2f}g$. Denoting by $\Omega ^{1}(M)$ the $C^{\infty }\left(
M\right) $-module of 1-forms on $M$ a (Riemannian) Weyl structure is a map $%
W:{\cal G}\rightarrow \Omega ^{1}\left( M\right) $ such that $W\left(
\overline{g}\right) =W\left( g\right) +2df$. In \cite{f:l} it is proved that
for a Weyl manifold $\left( M,{\cal G},W\right) $ there exists a unique
torsion-free linear connection $\nabla $ on $M$ such that for every $g\in
{\cal G}$:
$$
\nabla g=W\left( g\right) \otimes g.\eqno\left( \ast \right)
$$
The parallel transport induced by $\nabla $ preserves the given conformal
class ${\cal G}$.

A natural extension of Riemannian metrics are the generalized Lagrange
metrics. These metrics, introduced by Radu Miron around 1980, are suitable
in geometrical approaches of general relativity and gauge theory cf. \cite
{m:a}.

In this paper we extend the Weyl structures and compatible connections $%
\left( \ast \right) $ in the generalized Lagrangian framework.

\section{Generalized Lagrange-Weyl manifolds}

Let $M$ be a smooth $n$-dimensional real manifold and $\pi :TM\rightarrow M$
the tangent bundle of $M$. A chart $x=\left( x^{i}\right) _{1\leq i\leq n}$
of $M$ lifts to a chart $\left( x,y\right) =\left( x^{i},y^{i}\right) $ on $%
TM$. A tensor field of $\left( r,s\right) $-type on $TM$ with law of change,
at a change of charts on $M$, exactly as a tensor field of $\left(
r,s\right) $-type on $M$, is called {\it d-tensor field of }$\left(
r,s\right) ${\it -type}.

{\bf Definition 1.1}(\cite{m:a}) A d-tensor field of $\left( 0,2\right) $%
-type $g=\left( g_{ij}\left( x,y\right) \right) $ on $TM$ is called a {\it %
generalized Lagrange metric} ({\it GL-metric}, on short) if:

(i) it is symmetric: $g_{ij}=g_{ji}$

(ii) it is non-degenerate: $\det \left( g_{ij}\right) \neq 0$

(iii) the quadratic form $g_{ij}\left( x,y\right) \xi ^{i}\xi
^{j}$ has constant signature, $\forall \xi =\left( \xi ^{i}\right)
\in {I\!\!R}^{n}$.

{\bf Definition 1.2} Two GL-metrics $g,\overline{g}$ are called {\it %
conformal equivalent} if there exists $f\in C^{\infty }\left( M\right) $
such that $\overline{g}=e^{2f}g$.

In the following let ${\cal G}$ be a conformal structure i.e. an equivalence
class of conformal equivalent GL-metrics. The main notion of this paper is:

{\bf Definition 1.3} A {\it generalized Lagrange-Weyl structure} is a map $W:%
{\cal G}\rightarrow \Omega ^{1}\left( M\right) $ such that for every $g,%
\overline{g}\in {\cal G}$:
$$
W\left( \overline{g}\right) =W\left( g\right) +2df.\eqno\left( 1.1\right)
$$
The triple $\left( M,{\cal G},W\right) $ will be called {\it generalized
Lagrange-Weyl manifold}.

\medskip

Let us point that, from $\left( 1.1\right) $, if for some $g\in {\cal G}$
the 1-form $W\left( g\right) $ is closed (or exact) then for every $%
\overline{g}\in {\cal G}$ the 1-form $W\left( \overline{g}\right) $ is
closed (or exact).

On $TM$ the map $u\in TM\rightarrow V_{u}TM:=\ker \pi _{*,u}$ defines an
integrable distribution denoted $V\left( TM\right) $ and called {\it the
vertical distribution}. Recall that a vector field $X=X^{i}\left( x\right)
\frac{\partial }{\partial x^{i}}\in {\cal X}\left( M\right) $ has a {\it %
vertical lift} $X^{v}\in V\left( TM\right) $ given by $X^{v}=X^{i}\frac{%
\partial }{\partial y^{i}}$.

Because ${\cal G}$ implies the tangent bundle geometry it seems naturally
the following definition: a linear connection $\nabla $ on $TM$ is
vertical-compatible with the generalized Lagrange-Weyl structure $\left( M,%
{\cal G},W\right) $ if there exists $g\in {\cal G}$ such that for every $%
X\in {\cal X}\left( M\right) $:
\[
\nabla _{X^{v}}g=W\left( g\right) \left( X\right) \cdot g.
\]
But this definition has a great inconvenience: the fact that $\nabla $ is
vertical-compatible with a representative of ${\cal G}$ does not involve the
vertical-compatibility with another representative of ${\cal G}$. Indeed,
using $\left( 1.1\right) $, we have:
\[
\nabla _{X^{v}}\overline{g}=\nabla _{X^{v}}\left( e^{2f}g\right)
=X^{v}\left( e^{2f}\right) \cdot g+e^{2f}\nabla _{X^{v}}g=
\]
\[
=0\cdot g+e^{2f}W\left( g\right) \left( X\right) \cdot g=W\left( g\right)
\left( X\right) \cdot \overline{g}\neq W\left( \overline{g}\right) \left(
X\right) \cdot \overline{g}.
\]
With this motivation we introduce the next notion, namely {\it
nonlinear connections}, well-known in the geometry of tangent
bundle.

\section{Compatibility with respect to a nonlinear connection}

{\bf Definition 2.1}(\cite{m:a}) A distribution $H$ on $TM$ supplementary to
the vertical distribution i.e. $TTM=H\oplus V\left( TM\right) $ is called a
{\it nonlinear connection}.

An adapted basis for $V\left( TM\right) $ is $\left( \frac{\partial }{%
\partial y^{i}}\right) $ and an adapted basis for $H$ has the form $\left(
\frac{\delta }{\delta x^{i}}:=\frac{\partial }{\partial x^{i}}-N_{i}^{j}%
\frac{\partial }{\partial y^{j}}\right) $. The functions $\left(
N_{i}^{j}\left( x,y\right) \right) $ are called {\it the coefficients of the
nonlinear connection} $H$. We obtain a new lift for vector fields; namely,
to $X=X^{i}\left( x\right) \frac{\partial }{\partial x^{i}}\in {\cal X}%
\left( M\right) $ we associate the {\it horizontal lift} $X^{h}=X^{i}\frac{%
\delta }{\delta x^{i}}\in H$.

The nonlinear connection $H$ yields a bundle denoted $H(TM)$ and called {\it %
horizontal}. The existence of a nonlinear connection is equivalent to the
reduction of the standard almost tangent structure of $TM$ to a $D(GL(n,R))$%
-structure cf \cite{i:s}, \cite{a:i}.

{\bf Definition 2.2} A $D(GL(n,R))$-connection on $TM$ is called {\it %
d-connection} (or {\it Finsler connection}).

A d-connection $\nabla $ preserves by parallelism both the vertical and
horizontal bundles. Hence, $\nabla $ has a pair of Christoffel coefficients $%
\left( F_{jk}^{i}\left( x,y\right) ,C_{jk}^{i}\left( x,y\right) \right) $
defined by relations:
\[
\left\{
\begin{array}{c}
\nabla _{\frac{\delta }{\delta x^{j}}}\frac{\delta }{\delta x^{k}}=F_{jk}^{i}%
\frac{\delta }{\delta x^{i}},\quad \nabla _{\frac{\delta }{\delta x^{j}}}%
\frac{\partial }{\partial y^{k}}=F_{jk}^{i}\frac{\partial }{\partial y^{i}}
\\
\nabla _{\frac{\partial }{\partial y^{j}}}\frac{\delta }{\delta x^{k}}%
=C_{jk}^{i}\frac{\delta }{\delta x^{i}},\quad \nabla _{\frac{\partial }{%
\partial y^{j}}}\frac{\partial }{\partial y^{k}}=C_{jk}^{i}\frac{\partial }{%
\partial y^{i}}
\end{array}
\right. .
\]

It follows that $\nabla $ yields two algorithms of covariant derivation on
d-tensor fields: a horizontal one, denoted $_{|}$, and a vertical one,
denoted $|$. For example, on the d-tensor field $g=\left( g_{ij}\left(
x,y\right) \right) $ of $\left( 0,2\right) $-type we have:
$$
\left\{
\begin{array}{c}
g_{jk|i}=\frac{\delta g_{jk}}{\delta x^{i}}-g_{ak}F_{ji}^{a}-g_{ja}F_{ki}^{a}
\\
g_{jk}|_{i}=\frac{\partial g_{jk}}{\partial y^{i}}%
-g_{ak}C_{ji}^{a}-g_{ja}C_{ki}^{a}
\end{array}
\right. .\eqno\left( 2.1\right)
$$

{\bf Definition 2.3} A d-connection is called:

(i) {\it horizontal} if all $C_{jk}^{i}=0$,

(ii) {\it horizontal symmetric} ({\it h-symmetric} on short) if $%
F_{jk}^{i}=F_{kj}^{i}$ for all indices $i,j,k$,

(iii) {\it total symmetric} if it is h-symmetric and vertical symmetric i.e.
$C_{jk}^{i}=C_{kj}^{i}$ for all $i,j,k$.

For example, if $g$ is a Riemannian metric then the Levi-Civita connection
is the unique d-connection horizontal and h-symmetric; in this case $%
F_{jk}^{i}$ does not depend of $y$ since they are the usual Christoffel
coefficients.

It is natural to consider:

{\bf Definition 2.4} If $\left( M,{\cal G},W\right) $ is a generalized
Lagrange-Weyl manifold then a d-connection $\nabla $ is called {\it %
compatible} if there exists $g\in {\cal G}$ such that for every $X\in {\cal X%
}\left( M\right) $:

$$
\nabla _{X^{h}}g=W\left( g\right) \left( X\right) \cdot g.\eqno\left(
2.1\right)
$$

An important result is:

{\bf Proposition 2.5} {\it If }$\left( 2.1\right) ${\it \ holds for a given }%
$g\in {\cal G}${\it \ then }$\nabla ${\it \ is compatible with the whole
class} ${\cal G}$.

{\bf Proof} From $\left( 1.1\right) $ we get:
\[
\nabla _{X^{h}}\overline{g}=\nabla _{X^{h}}\left( e^{2\sigma }\circ \pi
\right) g=X^{h}\left( e^{2\sigma }\circ \pi \right) \cdot g+e^{2\sigma
}\nabla _{X^{h}}g=
\]
\[
=2d\sigma \left( X\right) e^{2\sigma }g+e^{2\sigma }W\left( g\right) \left(
X\right) g=
\]
\[
=e^{2\sigma }g\left( 2d\sigma +W\left( g\right) \right) \left(
X\right) =W\left( \overline{g}\right) \left( X\right) \cdot
\overline{g}.
\]

The pair $\left( g,H\right) $ yields four remarkable d-connections (\cite
{a:n}): Cartan, Berwald, Chern-Rund and Hashiguchi. For our aim, the
Chern-Rund connection, denoted $\nabla ^{CR}$, is more convenient because it
satisfies (\cite{a:n}):

I) is horizontal-metrical: $\nabla _{X^{h}}^{CR}g=0$ for every $X\in {\cal X}%
\left( M\right) $

II) is total symmetric.

The main result of this paper is:

{\bf Theorem 2.6} {\it For every generalized Lagrange-Weyl manifold }$\left(
M,{\cal G},W\right) ${\it \ there exists an unique compatible d-connection
which is horizontal and h-symmetric}.

{\bf Proof} Let $g\in {\cal G}$ and the associated $\nabla ^{CR}$. For $%
X,Y\in {\cal X}\left( M\right) $ let us define $\nabla _{X^{v}}Y^{v}=0$ and
:
$$
\nabla _{X^{h}}Y^{h}:=\nabla _{X^{h}}^{CR}Y^{h}-\frac{1}{2}W\left( g\right)
\left( X\right) \cdot Y^{h}-\frac{1}{2}W\left( g\right) (Y)\cdot X^{h}+\frac{%
1}{2}g\left( X^{h},Y^{h}\right) \cdot B\eqno\left( 2.2\right)
$$
where:
\[
g\left( X^{h},Y^{h}\right) =g_{ij}X^{i}Y^{j},\quad X=X^{i}\left(
x\right) \frac{\partial }{\partial x^{i}},Y=Y^{j}\left( x\right) \frac{%
\partial }{\partial x^{j}}
\]
and $B\in {\cal X}\left( TM\right) $ is $B=\left( B^{i}\right) $:
\[
B^{i}=g^{ij}wj,\quad W\left( g\right) =w_{i}dx^{i}.
\]
Here $\left( g^{ij}\right) $ is the inverse of $\left( g_{ij}\right) $.
Then:
\[
\nabla _{X^{h}}g=\nabla _{X^{h}}^{CR}g+W\left( g\right) \left(
X\right) \cdot g\stackrel{I)}{=}W\left( g\right) \left( X\right)
\]
i.e. $\nabla $ is horizontal-compatible with $g$. Applying the
previous result we have the conclusion.

\medskip

The non-null coefficients of $\nabla $ are:

\[
F_{jk}^{i}=\stackrel{CR}{F^{i}}_{jk}-\delta _{j}^{i}w_{k}-\delta
_{k}^{i}w_{j}+g_{jk}w^{i}
\]
where $\left( \stackrel{CR}{F}\right) $ are the horizontal Christoffel
coefficients of $\nabla ^{CR}$:

\[
\stackrel{CR}{F^{i}}_{jk}=g^{ia}\left( \frac{\delta g_{ak}}{\delta x^{j}}+%
\frac{\delta g_{ja}}{\delta x^{k}}-\frac{\delta g_{jk}}{\delta x^{a}}\right)
\]
and, as usual, $\left( w^{i}\right) $ is the $g$-contravariant version of $%
W\left( g\right) $ i.e. $w^{i}=g^{ia}w_{a}$.

\vspace{.3cm}

\noindent Faculty of Mathematics \newline University "Al. I. Cuza"
\newline Ia\c si, 700506 \newline Romania \newline email:
mcrasm@uaic.ro


\begin{thebibliography}{9}
\bibitem{a:i}  Aikou, T., Ichijyo, Y., {\it Finsler-Weyl structures and
conformal flatness}, Rep. Fac. Sci., Kagoshima Univ., 23(1990), 101-109.

\bibitem{a:n}  Anastasiei, M., {\it Finsler connections in generalized
Lagrange spaces}, Balkan J. Geom. Appl., 1(1996), no. 1, 1-9.

\bibitem{f:l}  Folland, G. B., {\it Weyl manifolds}, J. Diff. Geom.,
4(1970), 145-153.

\bibitem{i:s}  Ichijyo, I., {\it The }$D(O(n))${\it -structures in tangent
bundles}, J. Math. Tokushima Univ., 23(1989), 7-21.

\bibitem{k:l}  Kozma, L., {\it On Finsler-Weyl manifolds and connections},
Rend. Circ. Mat. Palermo (2), Suppl. No. 43(1996), 173-179.

\bibitem{m:a}  Miron, R., Anastasiei, M., {\it The Geometry of Lagrange
Spaces: Theory and Applications}, Kluwer Academic Publishers, FTPH, vol. 59,
1994.

\bibitem{w:h}  Weyl, H., {\it Space, time, matter}, Dover Publ., 1952.
\end{thebibliography}
\end{document}